\pdfoutput=1
\RequirePackage{ifpdf}
\ifpdf %We are running pdfTeX in pdf mode
\documentclass[pdftex]{sigma}
\else
\documentclass{sigma}
\fi

\numberwithin{equation}{section}
\newtheorem{Theorem}{Theorem}[section]

\newtheorem{Lemma}[Theorem]{Lemma}
\newtheorem{Proposition}[Theorem]{Proposition}
{\theoremstyle{definition}
\newtheorem{Definition}[Theorem]{Definition}

\newtheorem{Example}[Theorem]{Example}
\newtheorem{Remark}[Theorem]{Remark}
\newtheorem*{notation}{Notation}
}

\usepackage{overpic}

\begin{document}

\newcommand{\arXivNumber}{1410.7630}

\allowdisplaybreaks

\renewcommand{\PaperNumber}{008}

\FirstPageHeading

\ShortArticleName{Ambiguities in a~Problem in Planar Geodesy}

\ArticleName{Ambiguities in a~Problem in Planar Geodesy}

\Author{Josef SCHICHO and Matteo GALLET}

\AuthorNameForHeading{J.~Schicho and M.~Gallet}

\Address{Research Institute for Symbolic Computation, Johannes Kepler University,\\
Altenberger Strasse 69, 4040 Linz, Austria}

\Email{\href{mailto:josef.schicho@risc.jku.at}{josef.schicho@risc.jku.at}, \href{mailto:mgallet@risc.jku.at}{mgallet@risc.jku.at}}

\ArticleDates{Received October 29, 2014, in f\/inal form January 16, 2015; Published online January 27, 2015}

\Abstract{This is a~study of a~problem in geodesy with methods from complex algebraic geometry: for a~f\/ixed number of
measure points and target points at unknown position in the Euclidean plane, we study the problem of determining their
relative position when the viewing angles between target points seen from measure points are known.
In particular, we determine all situations in which there is more than one solution.}


\Classification{14N05; 14N15}

\section{Introduction}

\begin{figure}[ht] \centering
\begin{overpic}
[width=0.3\textwidth]{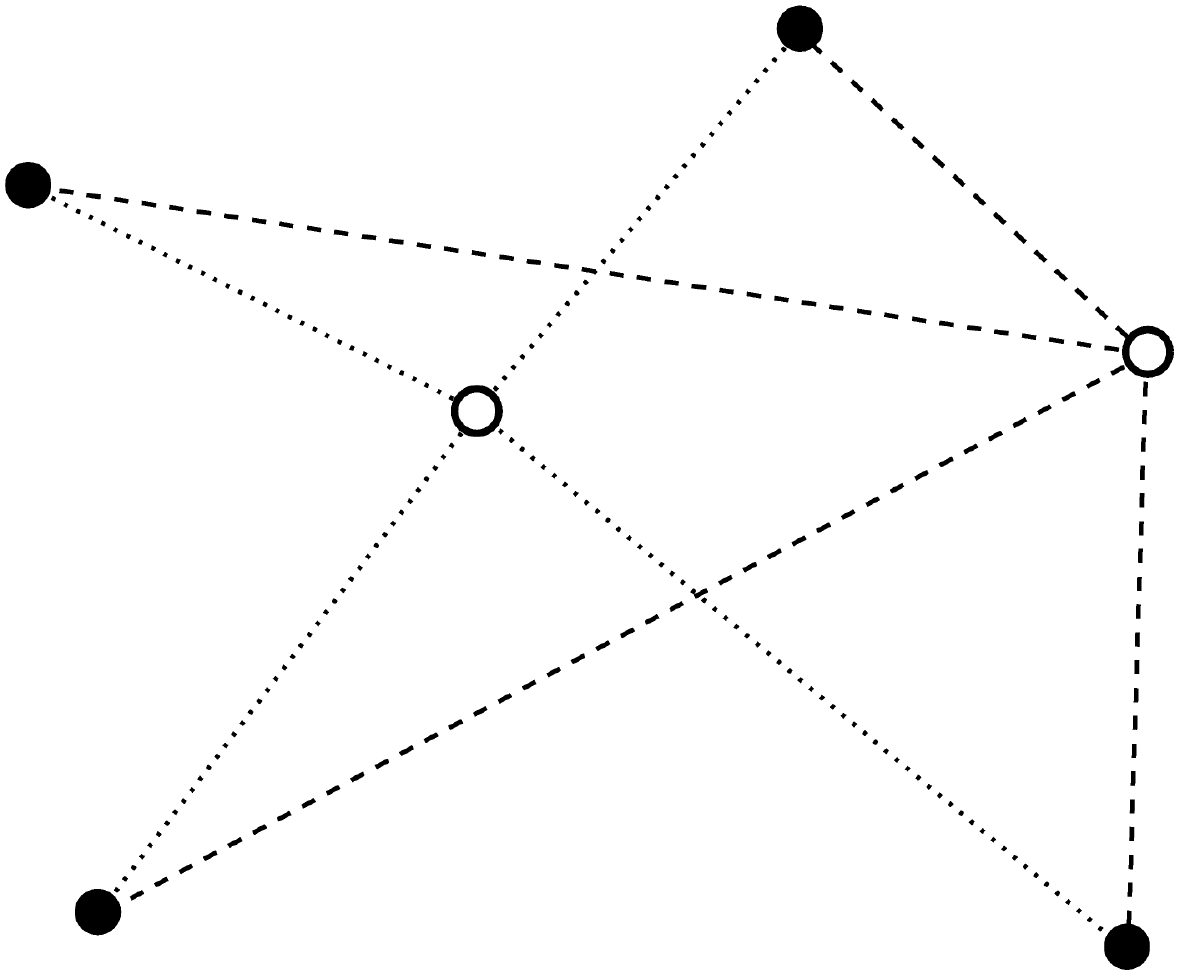}
\begin{small}
\put(0,0){$p_1$} \put(-8,66){$p_2$} \put(72,82){$p_3$} \put(100,0){$p_4$} \put(39,41){$q_1$} \put(102,51){$q_2$}
\end{small}
\end{overpic}
\caption{An example with $t = 4$ and $m = 2$, namely with four target points $p_1, \dotsc, p_4$ and two measure points
$q_1$, $q_2$.}
\label{figure:points_cameras}
\end{figure}
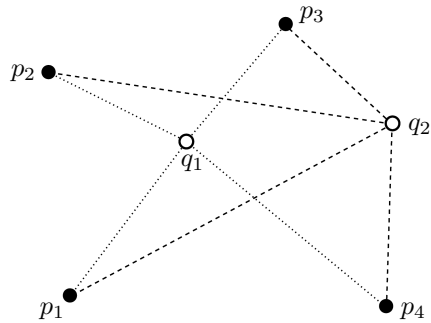

Let $t>0$ and $m>0$ be integers.
We consider the problem of identifying the relative position of $t+m$ unknown points $p_1,\dotsc,p_t$, $q_1,\dotsc,q_m$
in the Euclidean plane from the angles $\measuredangle_{q_j}(p_i,p_k)$, for $1\le i < k\le t$ and $1 \le j \le m$.
We call $p_1,\dotsc,p_t$ the {\em target points} and $q_1,\dotsc,q_m$ the {\em measure points} (see
Fig.~\ref{figure:points_cameras}).

Apparently, all angles can be computed from the angles $\measuredangle_{q_j}(p_1,p_{k})$ for $2 \le k \le t$ and $1\le
j\le m$, since $\measuredangle_{q_j}(p_i,p_{k}) = - \measuredangle_{q_j}(p_1,p_{i}) + \measuredangle_{q_j}(p_1,p_{k})$.
So the dimension of given data is $m(t-1)$, and the dimension of the unknown data is $2(m+t-2)$; in order that the
problem is not undetermined, we should have
\begin{gather*}
0\le m(t-1)-2(m+t-2)=(m-2)(t-3)-2.
\end{gather*}
In the two cases where equality holds, there are in general two solutions, namely there exist two $(m+t)$-tuples of
points in the plane both f\/itting with the data about the angles.

\looseness=-1
To explain the situation for $t = m = 4$, we introduce the {\em profile} of a~f\/ixed sequence $\vec{p}$ of target points,
as the two-dimensional set of all possible measurement results from various measure points.
Then we decompose our problem into two steps.
First, we interpolate the prof\/ile surface from the given angles.
Once the prof\/ile is determined, additional measure points (and the corresponding measured angles) do not reveal any more
information on the position of the target points.
Then we determine the target points from the prof\/ile.
Remarkably, it turns out that for target sequences~$\vec{p}$ with $t = 4$, there is in general a~second sequence~$\vec{p}\,{}'$ with the same prof\/ile.

Similarly, we treat the case $m = 3$ and $t = 5$ by introducing the {\em co-profile} of a~sequen\-ce~$\vec{q}$ of measure
points, namely the two-dimensional set of measurements of target points from~$\vec{q}$; however, the situation is not
entirely symmetric to the case of the prof\/ile because here we need in addition one target point in order to def\/ine the
co-prof\/ile.
In this case, additional target points (and the corresponding measured angles) do not reveal any more information on the
position of the measure points and the prescribed target point.

Section~\ref{sec:ad} is devoted to the introduction of prof\/ile and co-prof\/ile for arbitrary~$t$ and~$m$.
In Section~\ref{sec:ipol}, we determine the situations where the prof\/ile or the co-prof\/ile cannot be uniquely
determined, also paying attention to non-general situations.
In Section~\ref{sec:point_identification} we identify the target points once the prof\/ile is known.
We determine the situations where this is not uniquely possible.
Once we know the target points, the determination of the measure points from the prof\/ile is not problematic.
The complete discussion is summarized in Theorem~\ref{thm:sup}.

Strictly speaking, the two solutions for the general case $m = t = 4$ cannot be distinguished algebraically, by means of
testing polynomial equalities.
But angle measurement also gives information corresponding to inequalities, namely we can detect on which of the two
rays of a~line a~certain target point lies.
Taking this additional information into account, one can distinguish the two solutions in many cases, but not in all of
them.
This aspect is discussed in more detail in Section~\ref{sec:twinmap}.

\looseness=-1
The majority of the results in this paper is not new.
The stated problem is equivalent to the reconstruction of a~set of points in the plane from images by a~one-dimensional
calibrated camera.
The paper~\cite{AstromOskarsson} (also contained, in expanded form, in~\cite{Oskarsson}) explains the two solutions for
the case $m=3$ by a~theorem of planar geometry.
There the authors also describe a~duality exchanging the role of measure points and target points, except for one
target point that has to keep its role as a~target point; using this duality, the existence of two solutions for $t=4$
is explained as well.
A~similar duality is also known for the uncalibrated case (see~\cite{Quan}), where three f\/ixed target points keep their
role instead of one, and also for two-dimensional pictures of points in 3-space (see~\cite{CarlssonWeinshall}).
The two solutions for the general case appear also in the one-dimensional uncalibrated case in~\cite{Quan}; indeed, the
calibrated case can be reduced to the uncalibrated case by considering the two cyclic points of the Euclidean plane as
additional target points.
The classif\/ication of exceptional cases with inf\/initely many solutions in Theorem~\ref{thm:sup} seems to be new, as well
as the discussion of ambiguities taking direction information into account (Theorem~\ref{thm:real}).
However, we think that our uniform discussion of the ambiguities has some value.
Also, the results are obtained by general theorems in algebraic geometry, without technical computations; a~single
exception is Remark~\ref{rem:twin}, which is not used in the remaining part of the paper.

\section{Algebraic description of given data}
\label{sec:ad}

\begin{notation}
Throughout the paper we will often use the adjective \emph{general} talking about points or measurements: as it is often
done in algebraic geometry, here we mean that a~property is \emph{general} if it holds for a~Zariski open set of a~space
parametrizing the objects we are interested in.
Hence when we say that a~property holds for ``general~$t$ target points~$\vec{p}$ and~$m$ measure points~$\vec{q}\,$'' we
mean that the set of pairs $(\vec{p}, \vec{q})$ for which the property does not hold is a~variety of dimension less
than~$2(m+t)$.
\end{notation}

\vspace{-1mm}

The aim of this section is to associate to any given tuple~$\vec{p}$ of target points a~projective variety, called the
\emph{profile} of~$\vec{p}$, and to any given tuple~$\vec{q}$ of measure points together with a~prescribed target
point~$p_1$ another projective variety, called the \emph{co-profile} of~$\vec{q}$ and~$p_1$.
These two varieties encode information about, respectively, all possible measurements of~$\vec{p}$ from any measure
points, and measurements with respect to~$\vec{q}$ and~$p_1$ of any possible target point.

\begin{notation}
For target points, we use real coordinates $(u,v)$, and we introduce complex coordinates $w= u + \mathrm{i} v$ and
$\overline{w} = u - \mathrm{i} v$.
For measure points, we use real coordinates $(x,y)$, and we introduce complex coordinates $z = x + \mathrm{i} y$ and
$\overline{z} = x - \mathrm{i} y$.
\end{notation}

The f\/irst aspect to clarify is how we model measurements, namely angles.
We decide to model them by complex numbers with modulus~$1$.
More precisely, if the measure point~$q$ has complex coordinate~$z$ and the target points~$p_1$ and~$p_2$ have complex
coordinates~$w_1$ and~$w_2$, one can check that the angle $\measuredangle_{q}(p_1,p_2)$ is given by the complex number
$\frac{(\overline{z} - \overline{w}_1)(z - w_2)}{\left| (\overline{z} - \overline{w}_1)(z - w_2) \right|}$.
We f\/ix now pairwise distinct target points $p_1,\dotsc,p_t$ (with real/complex coordinates as above).
Hence we can def\/ine a~map $\mathbb R^2 \setminus \{p_1, \dotsc, p_t\} \longrightarrow \mathbb C^{t-1}$ sending
\begin{gather*}
z = x + \mathrm{i} y
\
\mapsto
\
\left(\frac{(\overline{z} - \overline{w}_1)(z - w_2)}{\left| (\overline{z} - \overline{w}_1)(z - w_2) \right|}, \dotsc,
\frac{(\overline{z} - \overline{w}_1)(z - w_{t})}{\left| (\overline{z} - \overline{w}_1)(z - w_{t}) \right|} \right).
\end{gather*}
This map is far from being algebraic.
In order to make it algebraic, we start by squaring each coordinate of the image (this corresponds to multiplying all
angles by a~factor of~$2$) and then homogenize with respect to a~new variable, which we place in f\/irst position, so that
the domain becomes a~projective space.
What we obtain has the following expression, which still involves complex conjugation, so it is not yet what we are
looking for:
\begin{gather*}
f_{\vec{p}}:  \ \mathbb R^2 \setminus \{p_1,\dotsc,p_t\}   \longrightarrow   \mathbb P_{\mathbb C}^{t-1},
\\
\hphantom{f_{\vec{p}}: \ }~z   \mapsto   \big(F_1(z): \ldots: F_t(z) \big),
\\
F_1(z)   =   (z-w_1)(\overline{z} -\overline{w}_2) \cdots (\overline{z}-\overline{w}_t),
\\
F_2(z)   =   (\overline{z}-\overline{w}_1)(z -w_2) \cdots (\overline{z}-\overline{w}_t),
\\
\cdots\cdots\cdots\cdots\cdots\cdots\cdots\cdots\cdots\cdots\cdots\cdots
\\
F_t(z)   =   (\overline{z}-\overline{w}_1)(\overline{z} -\overline{w}_2)\cdots (z-w_t).
\end{gather*}

We proceed towards our goal to obtain an algebraic map, and in particular we want to make~$f_{\vec{p}}$ into a~rational
map between projective varieties.
To do so we choose the following injection\footnote{The choice of $\mathbb P^1_{\mathbb C} \times \mathbb P^1_{\mathbb
C}$, instead of, for example, $\mathbb P^2_{\mathbb C}$, will be justif\/ied \emph{a posteriori} by the fact that in this
way we will be able to prove some properties of the map~$f_{\vec{p}}$ in an easier way.} of~$\mathbb R^2$ into $\mathbb
P^1_{\mathbb C} \times \mathbb P^1_{\mathbb C}$:
\begin{gather*}
\mathbb R^2  \longrightarrow  \mathbb P^1_{\mathbb C} \times \mathbb P^1_{\mathbb C},
\\
(x,y)  \mapsto  \underbrace{(x + \mathrm{i} y: 1) \times (x - \mathrm{i} y: 1)}_{= (z: 1) \times (\overline{z}:1)}.
\end{gather*}
Then we can extend the previously def\/ined map~$f_{\vec{p}}$ to a~map from $\mathbb P^1_{\mathbb C} \times \mathbb
P^1_{\mathbb C}$, which we still denote by $f_{\vec{p}}$.
If we take coordinates $(\alpha_1:\beta_1) \times (\alpha_2:\beta_2)$ on $\mathbb P^1_{\mathbb C} \times \mathbb
P^1_{\mathbb C}$, the new map is given by
\begin{gather*}
f_{\vec{p}}:  \ \mathbb P^1_{\mathbb C} \times \mathbb P^1_{\mathbb C}   \dashrightarrow   \mathbb P_{\mathbb C}^{t-1},
\\
\hphantom{f_{\vec{p}}: \ }~(\alpha_1:\beta_1) \times (\alpha_2:\beta_2)   \mapsto   \big(G_1(\alpha, \beta): \ldots:G_t(\alpha, \beta)\big).
\end{gather*}
The components~$G_j$ of~$f_{\vec{p}}$ are obtained from the previously def\/ined polynomials~$F_j$ by substituting all
factors $(z - w_i)$ with $(\alpha_1 - \beta_1 w_i)$ and all factors $(\overline{z} - \overline{w}_i)$ with $(\alpha_2 -
\beta_2 \overline{w}_i)$.
Thus we see that this is a~rational map between complex projective varieties.

\begin{Remark}
Note that, since during the construction of~$f_{\vec{p}}$ we take the squares of the involved quantities (this
corresponding to multiplication by~$2$ of the angles), we have that~$f_{\vec{p}}(q)$ determines the measured angles only
modulo~$\pi$.
This is a~consequence of the use of complex numbers: since complex numbers are not ordered, the ``complex-valued'' part
of our angle measurement device is unable to tell two opposite directions of the same line apart.
\end{Remark}

\begin{Remark}
If we remove the last target point and we consider $\vec{p^{*}} = (p_1, \dotsc, p_{t-1})$, then the map
$f_{\vec{p^{*}}}$ is equal to $f_{\vec{p}}$ composed with the projection $\mathbb P_{\mathbb C}^{t-1}\to \mathbb
P_{\mathbb C}^{t-2}$ given by $(s_1: \ldots: s_t) \mapsto (s_1: \ldots: s_{t-1})$.
\end{Remark}

\begin{Definition}
%\label{defn:profile}
The Zariski closure of the image of the map~$f_{\vec{p}}$ is an algebraic surface $S_{\vec{p}}\subset \mathbb P_{\mathbb
C}^{t-1}$, which we call the {\em profile} of $\vec{p}$.
\end{Definition}

\begin{Proposition}
\label{prop:measure_birational}
If $t \ge 3$, then the measurement map~$f_{\vec{p}}$ is a~birational map from the plane to the profile of~$\vec{p}$, for
any point sequence~$\vec{p}$ of pairwise distinct points.
\end{Proposition}

\begin{proof}
First, assume that $t = 3$.
We recall a~well-known (for surveying or land navigation purposes; see for example~\cite{Mooers}) geometric construction
of 	a left inverse of the measurement map: for each $q \in \mathbb R^2$, the oriented angle $\measuredangle_q(p_1,p_2)$
modulo~$\pi$ determines a~circle on which~$q$ lies.
Similarly, the oriented angle $\measuredangle_q(p_1,p_3)$ modulo~$\pi$ determines another circle.
The two circles intersect in~$p_1$ and in a~	second point, which must be~$q$.
It is straightforward to express this geometric construction as a~rational map in the angles.

For $t > 3$, a~left inverse can be given as the projection to the f\/irst three coordinates followed by the inverse for
the case $t = 3$ above.
\end{proof}

\begin{notation}
The inverse of the measurement map is called the {\em resection map}.
\end{notation}

\begin{Remark}
%\label{rem:exceptional}
From the geometric construction above, it is clear that the measurement map is injective outside the circle through any
triple of points $p_i$, $p_j$, $p_k$ (or the line through $p_i$, $p_j$, $p_k$ if the points are collinear).
\end{Remark}

Here is a~description of the prof\/ile.

\begin{Lemma}
\label{lem:degree_profile}
If $t \ge 3$, then the profile surface $S_{\vec{p}} \subset \mathbb P_{\mathbb C}^{t-1}$ has degree $t-2$ and passes
through the $t+1$ points $(1:\ldots:0), \ldots, (0:\ldots:1), (1:\ldots:1)$.
Its ideal is generated by the $2 \times 2$ minors of a~$2 \times (t-2)$ matrix whose entries are linear forms.
\end{Lemma}

\begin{proof}
The functions $G_1, \dotsc, G_t$ have bidegree $(1, t-1)$.
The base locus of $f_{\vec{p}}$ consists of~$t$ simple base points corresponding to $p_1, \dotsc, p_t$, hence the degree
of the image is $2(t-1)-t = t-2$.
For an explanation regarding the previous formula we refer to~\cite[Appendix~A]{Cox}\footnote{Here the author
clarif\/ies~\cite[equation~(5.1)]{Cox} using arguments from intersection theory.
The only dif\/ference between our formula and~\cite[equation~(5.1)]{Cox} is that in the latter one considers homogeneous
polynomials of degree~$n$, and this justif\/ies the term $n^2$ appearing there, while in our formula we have bihomogeneous
polynomials of bidegree $(1, t-1)$, yielding the term~${2(t-1)}$.
The term~$t$ in our formula, given by the~$t$ simple base points, corresponds to the term $\sum\limits_{p \in Z}
e(I_{Z,p}, \mathcal{O}_{\mathbb P^2,p})$ in~\cite[equation~(5.1)]{Cox}.}.
For $r \in \{1, \dotsc, t\}$, the image of the line $\alpha_2 - \beta_2 \overline{w}_r$ is the point
$(0:\ldots:1:\ldots:0)$ (with~$1$ at position~$r$); the point $(1:\ldots:1)$ appears as the image of the point at
inf\/inity $(1:0) \times (1:0)$.

Because $\deg(S_{\vec{p}}) = \mathrm{codim}(S_{\vec{p}}) + 1$, we have a~surface of minimal degree.
Such a~surface is determinantal, and the generators of its ideal are well-known (see for example~\cite{EisenbudHarris}).
\end{proof}

We come to the def\/inition of the \emph{co-profile}: for this purpose let us f\/ix a~tuple $\vec{q} = (q_1, \dotsc, q_m)$
of measure points and a~target point~$p_1$.
We consider the map $\mathbb R^2 \longrightarrow \mathbb C^{m}$ associating to each target point $p \in \mathbb R^2$
the~$m$-tuple of angles $\big(\measuredangle_{q_1}(p_1,p), \dotsc, \measuredangle_{q_m}(p_1,p)\big)$.
Performing analogous homogeneization and operations as before we obtain a~rational map:
\begin{gather*}
f'_{\vec{q},p_1}: \ \mathbb P^1_{\mathbb C} \times \mathbb P^1_{\mathbb C}  \dashrightarrow  \mathbb P^{m}_{\mathbb C},
\\
\hphantom{f'_{\vec{q},p_1}: \ }~(\alpha_1:\beta_1) \times (\alpha_2:\beta_2)   \mapsto   \big(G'_0(\alpha, \beta): \ldots:G'_m(\alpha, \beta)\big)
\end{gather*}

\begin{Definition}
%\label{defn:coprofile}
The Zariski closure of the image of the map~$f'_{\vec{q},p_1}$ is an algebraic surface $S'_{\vec{q},p_1}\subset \mathbb
P_{\mathbb C}^{m}$, which we call the {\em co-profile} of $\vec{q}$ and $p_1$.
\end{Definition}

\begin{Lemma}
%\label{lem:degree_coprofile}
If $m \ge 2$, then the co-profile map~$f'_{\vec{q},p_1}$ is birational, and the co-profile surface $S'_{\vec{q},p_1}
\subset \mathbb P_{\mathbb C}^{m}$ has degree $m-1$ and passes through the $m+1$ points $(1:\ldots:0), \ldots,
(0:\ldots:1), (1:\ldots:1)$.
Its ideal is generated by the $2 \times 2$ minors of a~$2 \times (m-1)$ matrix whose entries are linear forms.
\end{Lemma}
\begin{proof}
For the inverse of~$f'_{\vec{q},p_1}$, there is again a~geometric construction: when we know
$\measuredangle_{q_1}(p_1,p)$ and $\measuredangle_{q_2}(p_1,p)$, then we can construct~$p$ by intersecting the two lines
determined by the angles above.

The proof of the statement on $S'_{\vec{q},p_1}$ is analogous to the proof of Lemma~\ref{lem:degree_profile}, with some
changes in the details: the functions $G_0', \dotsc, G_m'$ have bidegree $(1,m)$, the base points correspond to $q_1,
\dotsc, q_m$ and the point at inf\/inity $(1:0) \times (1:0)$; the images of the exceptional lines are the points
$(1:\ldots:0), \ldots, (0:\ldots:1)$, and the point $(1:\ldots:1)$ appears as the image of~$p_1$.
\end{proof}

\section{An interpolation problem}
\label{sec:ipol}

\looseness=-1
As remarked in the Introduction, we are going to solve our problem in two steps.
The f\/irst step, discussed in this section, is to compute the prof\/ile of a~tuple~$\vec{p}$ from a~f\/ixed number~$m$ of
measurements.
The task is to f\/ind a~surface of known degree, and we know $m+t+1$ points on the surface, namely~$m$ points from
measurements and~$t+1$ points from Lemma~\ref{lem:degree_profile}.
We call this step \emph{profile interpolation}.
Once the prof\/ile is computed, we know all possible measurements up to~$\pi$, hence additional values obtained with the
``complex angle measurement device'', namely the map~$f_{\vec{p}}$, do not provide more information.
Analogously, one can try to compute the co-prof\/ile of~$\vec{q}$ and~$p_1$ from other~$t-1$ target points.
We will see that the solution of this problem is not always unique, namely there may exist more than one
prof\/ile/co-prof\/ile matching the given data.

If we are given~$t$ target points~$\vec{p}$ and~$m$ measure points~$\vec{q}$, then we can form an $m \times (t-1)$
matrix $M_{\vec{p}, \vec{q}}$, called the \emph{double angle matrix}:
\begin{gather*}
M_{\vec{p},\vec{q}}
=
\left(
\begin{matrix}
\big(\measuredangle_{q_1}(p_1,p_2)\big)^2 & \dots & \big(\measuredangle_{q_1}(p_1,p_t)\big)^2
\\
\vdots & & \vdots
\\
\big(\measuredangle_{q_m}(p_1,p_2)\big)^2 & \dots & \big(\measuredangle_{q_m}(p_1,p_t)\big)^2
\end{matrix}
\right).
\end{gather*}
In this way the rows of~$M_{\vec{p},\vec{q}}$ are the input for prof\/ile identif\/ication, while the columns
of~$M_{\vec{p},\vec{q}}$ are the input for co-prof\/ile identif\/ication.

We start with the case $t = 4$.
Here the prof\/ile is a~quadric surface in~$\mathbb P_{\mathbb C}^3$.

\begin{Proposition}
\label{prop:t4int}
If $t = 4$, then the profile can be computed from~$4$ measurements, in general.
\end{Proposition}

\begin{proof}
The linear space of quadratic forms in~$4$ variables is $10$-dimensional, and each point gives a~linear condition on it.
For general measurements, the~$4$ points in~$S_{\vec{p}}$ from the measurements and the~$5$ points from
Lemma~\ref{lem:degree_profile} give linear independent conditions, and so there is a~one-dimensional solution space
which determines the prof\/ile uniquely.
\end{proof}

\begin{Remark}
\label{rem:int}
If $t > 4$, then it is also possible to compute the prof\/ile from~$4$ general measurements.
For instance, if $\vec{p}=(p_1,\dotsc,p_5)$, then we can compute the prof\/iles~$S_1$ and~$S_2$ of $(p_1,p_2,p_3,p_4)$ and
of $(p_1,p_2,p_3,p_5)$, respectively.
Then the prof\/ile of~$\vec{p}$ is the closure of the set
\begin{gather*}
\big\{(x_0:x_1:x_2:x_3:x_4)\in \mathbb P_{\mathbb C}^4 \,\big|\,  (x_0:x_1:x_2:x_3)\in S_1~\text{and}~(x_0:x_1:x_2:x_4)\in S_2 \big\}.
\end{gather*}
\end{Remark}

As the parameter counting argument in the Introduction suggests, we can remove one measure point if $t \ge 5$ and still
have a~f\/inite number of solutions.
Here is the precise statement.

\begin{Proposition}
\label{prop:t5int}
For $t=5$ and $m=3$, there are in general two solutions for the profile interpolation problem.
\end{Proposition}

\begin{proof}
In this case the prof\/ile is a~rational cubic surface in~$\mathbb P^4_{\mathbb C}$.
By~\cite[Example~3.4]{EisenbudPopescu}, there are exactly two such cubic surfaces passing through $9$ points in general
position in~$\mathbb P_{\mathbb C}^4$.
Interestingly, the proof there uses Gale duality, which seems to be tightly related to the correspondence between
the~$m$ points on the prof\/ile and the $t-1$ points on the co-prof\/ile given by a~double angle
matrix~$M_{\vec{p},\vec{q}}$: in both cases, in fact, the coordinates of the two sets of points are given by the rows of
a~matrix and its transpose.
\end{proof}

Intuitively, increasing~$t$ should increase the amount of information, so we may think that for $t > 5$ and $m = 3$, we
have only a~single solution for the prof\/ile interpolation problem.
We will see later that this is not the case.
What is true is that considering also the co-prof\/ile increases the amount of insight.

\begin{Proposition}
\label{prop:m3int}
If $m=3$ and $t\ge 5$, then the co-profile can be uniquely computed, in general.
\end{Proposition}

\begin{proof}
The statement is analogous to Proposition~\ref{prop:t4int}: again, we have to interpolate~$9$ or more points via
a~quadric surface in~$\mathbb P_{\mathbb C}^3$.
\end{proof}

\begin{Proposition}
\label{prop:m4int}
For $t=4$ and $m=4$, there are in general two solutions for the co-profile interpolation problem.
\end{Proposition}

\begin{proof}
Here the situation is as in Proposition~\ref{prop:t5int}: we have to interpolate~$9$ points using cubic surfaces
in~$\mathbb P^4_{\mathbb C}$.
\end{proof}

\begin{Proposition}
%\label{prop:surprise}
For $m = 3$ and $t > 5$, there are in general two solutions for the profile interpolation problem.
Dually, for $m > 4$ and $t = 4$, there are in general two solutions for the co-profile interpolation problem.
\end{Proposition}

\begin{proof}
If $m = 3$ and $t = 5$, then we know two possible candidates for the prof\/ile, giving two possible candidates for the
position of measure/target points; however the two co-prof\/iles must be equal by Proposition~\ref{prop:m3int}.
But if we know the co-prof\/ile, then any measurements from new target points does not give any additional information
about the position of the old target points.
Therefore we will never be able to tell the true prof\/ile, no matter how many new target points we introduce.

The proof of the dual statement is analogous.
\end{proof}

So far we completed the description of the prof\/ile/co-prof\/ile interpolation in the general case, and we summarize the
results in Theorem~\ref{thm:summarize_profile}.

\begin{Theorem}
\label{thm:summarize_profile}
Let~$\vec{p}$ and~$\vec{q}$ be general~$t$ and~$m$-tuples of points in $\mathbb R^2$.
Then the following tables summarize the number of possible profiles and co-profiles compatible with the double angle
matrix~$M_{\vec{p},\vec{q}}$:

\begin{table}[ht]\centering
\begin{minipage}{0.45\textwidth} \centering \caption{Number of prof\/iles f\/itting, in general, a~given
double angle matrix.}\vspace{1mm}
\begin{tabular}{c|ccc}
${m} \ \backslash \  {t}$ & 3 & 4 & 5 {or more}
\\
\hline
3 & 1 &~$\infty$ & 2
\\
4 {or more} & 1 & 1 & 1
\\
\end{tabular}
\label{table:profiles}
\end{minipage}
\qquad
\begin{minipage}{0.45\textwidth} \centering \caption{Number of co-prof\/iles f\/itting, in ge\-ne\-ral,
a~given double angle matrix.}\vspace{1mm}
\begin{tabular}{c|ccc}
${m} \ \backslash \ {t}$ & 3 & 4 & 5 {or more}
\\
\hline
3 &~$\infty$ &~$\infty$ & 1
\\
4 {or more} &~$\infty$ & 2 & 1
\\
\end{tabular}
\label{table:coprofiles}
\end{minipage}
\end{table}
\end{Theorem}

\looseness=-1
In Theorem~\ref{thm:cyccub} we analyze those cases for which for a~general choice of the points~$\vec{p}$ and~$\vec{q}$
there is a~unique prof\/ile or co-prof\/ile, but for some special choices there may be more than one prof\/ile/co-prof\/ile, and
we identify how these special choices look like.
Considering Tables~\ref{table:profiles} and~\ref{table:coprofiles} we see that, in general, we have uniqueness for
the prof\/ile if $m = t = 3$ or if $m \geq 4$ and $t \geq 4$, and for the co-prof\/ile if $m \geq 3$ and $t \geq 5$.
Notice that if $t = 3$, then the prof\/ile is equal to $\mathbb P^2_{\mathbb C}$, and so in this case it can never happen,
even in special situations, that the prof\/ile is not unique.

\begin{notation}
A~\emph{cyclic cubic curve} is a~plane cubic curve passing through the two cyclic points at inf\/inity, with
homogeneous coordinates $(0:1:\pm\mathrm{i})$.
In the embedding of $\mathbb R^2$ into $\mathbb P^1_{\mathbb C} \times \mathbb P^1_{\mathbb C}$ that we use, the
equation of a~cyclic curve has bidegree $(2,2)$ and passes through the point $(1:0) \times (1:0)$.

We allow cyclic cubics to be reducible.
There are two reducible cases: a~line and a~circle, or a~conic and the line at inf\/inity.
\end{notation}

We start treating the base cases.

\begin{Lemma}
\label{lem:base_cases}
Suppose that $t = 4$ and $m \geq 4$ and that the profile is not uniquely determined by the double angle matrix
$M_{\vec{p},\vec{q}}$, or that $t \geq 5$ and $m = 3$ and the co-profile is not uniquely determined by the double angle
matrix $M_{\vec{p},\vec{q}}$.
Then there exists a~cyclic cubic curve that contains all target points and all measure points, or a~line that contains
all target points.
\end{Lemma}
\begin{proof}
First, assume $t=4$ and $m \geq 4$, so that the prof\/ile surface~$S$ is a~quadric surface in~$\mathbb{P}^3_{\mathbb{C}}$.
Moreover we assume that the prof\/ile is not unique.
This is equivalent to the fact that the $m+5$ points $f_{\vec{p}}(q_1), \dotsc, f_{\vec{p}}(q_m)$ and
$(1:\ldots:0),\dotsc,(0:\ldots:1)$ and $(1:\ldots:1)$ are contained in the intersection with a~second quadric~$Q$.
The pullback under $f_{\vec{p}}$ of any plane section of~$S$ is a~divisor of $\mathbb P_{\mathbb C}^1 \times \mathbb
P_{\mathbb C}^1$ of bidegree $(1,3)$ passing through $p_1, \dotsc, p_4$, therefore the pullback of any quadric section
is a~divisor of bidegree $(2,6)$ passing with multiplicity at least~$2$ through $p_1, \dotsc, p_4$.
In addition, by hypothesis~$Q$ passes through the point $(1:\ldots:0)$, which is the image of an exceptional
divisor~$E_1$, namely the unique curve of bidegree $(0,1)$ passing through~$p_1$.
Hence the pullback of~$Q$ contains this divisor~$E_1$.
Similarly, the pullback contains the other~$3$ exceptional divisors $E_2$, $E_3$ and $E_4$.
If we remove these components, it remains a~curve of bidegree $(2,2)$ passing through $p_1,\dotsc,p_4$.
It must contain the preimage of $(1:\ldots:1)$, which is the point $(1:0) \times (1:0)$, and all measure points
$q_1,\dotsc,q_m$.
So we have found a~cyclic cubic curve as stated in the thesis.

The proof for the dual assertion assuming non-uniqueness of the co-prof\/ile is similar.
\end{proof}

\begin{Example}
Consider the case $t = m = 4$ where all~$4$ target points $p_1, \dotsc, p_4$ lie on a~circle and the $2$ measure points  %, where
$q_1$ and $q_2$ lie on the same circle, as in Fig.~\ref{figure:profile_not_unique}.
In this case all~$8$ points lie on a~(reducible) cyclic cubic, given by the union of the circle and a~line.
\begin{figure}[ht] \centering
\begin{overpic}
[width=0.3\textwidth]{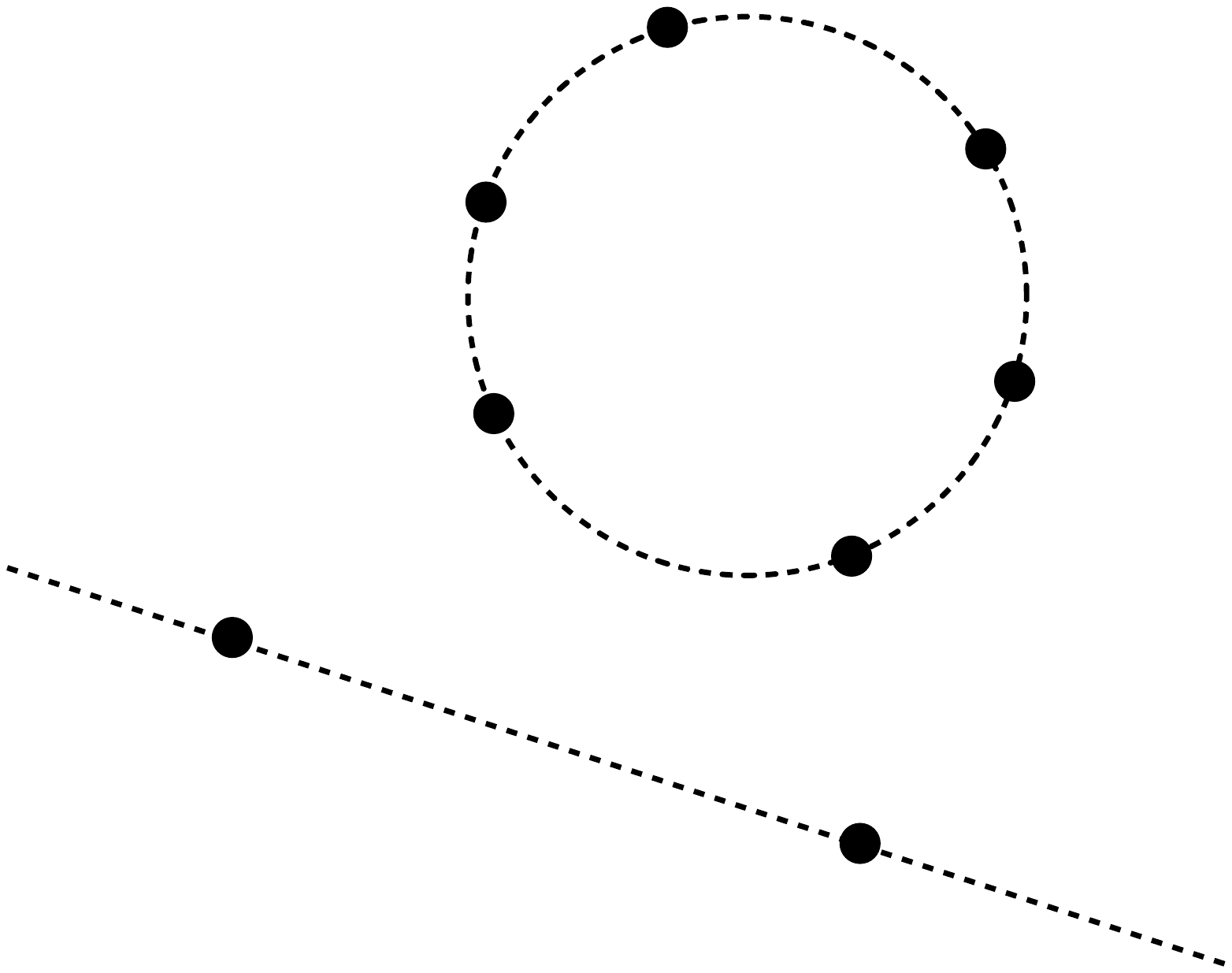}
\begin{small}
\put(31,59){$p_1$} \put(50,70){$p_2$} \put(74,61){$p_3$} \put(84.5,44.5){$p_4$} \put(34,40){$q_1$} \put(70,28){$q_2$}
\put(16,20){$q_3$} \put(64,4){$q_4$}
\end{small}
\end{overpic}
\caption{An example of a~conf\/iguration of target and measure points for which the prof\/ile is not unique: the four target
points lie on a~circle, two measure points lie on the same circle, and the other two measure points lie on a~line.}
\label{figure:profile_not_unique}
\end{figure}

Since all $p_1, \dotsc, p_4,q_1$ and $q_2$ lie on the same circle, we have the equalities
\begin{gather*}
\measuredangle_{q_1}(p_1,p_k)
=
\measuredangle_{q_2}(p_1,p_k)
\qquad
\text{for all}
\ \
k \in \{2, 3, 4 \},
\end{gather*}
which imply that the conditions imposed by $q_1$ and $q_2$ to the linear space of quadratic forms in $4$ variables as in
the proof of Proposition~\ref{prop:t4int} are not independent.
This means that there exists at least a~two-dimensional linear space of such forms, this ensuring the existence of
inf\/initely many possible prof\/iles for the pair $(\vec{p}, \vec{q})$.
\end{Example}

It would be nice now to use an inductive argument (on the number of target or measure points) to establish the same
result as in Lemma~\ref{lem:base_cases} for all remaining cases.
Unfortunately, in each induction step we encounter new possibilities, making a~uniform description dif\/f\/icult to achieve.
Because of this, we prefer to strengthen the hypotheses of our theorem in order to get rid of these ``spurious'' cases.

\begin{Theorem}
\label{thm:cyccub}
Assume that, for given sequences~$\vec{p}$ of target points and~$\vec{q}$ of measure points, either $m \ge 4$ and the
profile is not uniquely determined by the double angle matrix $M_{\vec{p},\vec{q}}$, or $t \ge 5$ and the co-profile is
not uniquely determined by the double angle matrix $M_{\vec{p},\vec{q}}$.
Assume furthermore that no $5$ target points and no~$4$ measure points are cocircular or collinear.
Then there exists a~cyclic cubic curve that contains all target points and all measure points.
\end{Theorem}

\begin{proof}
Let us start analyzing the case when $t=5$ and $m \geq 4$.
We try to reduce to Lemma~\ref{lem:base_cases}.
By hypothesis we know that the prof\/ile is not unique.
If now we consider any subset $\vec{p}\,{}'$ of four points out of the f\/ive points of $\vec{p}$, we have that the prof\/ile
of~$\vec{p}\,{}'$ and~$\vec{q}$ may or may not be unique.
Notice that if we could f\/ind two four tuples~$\vec{p}\,{}'$ and~$\vec{p}\,{}''$ for which the prof\/ile is unique, then~by
Remark~\ref{rem:int} we would be able to reconstruct the prof\/ile of $\vec{p}$ and~$\vec{q}$ uniquely, and this is
against the assumption.
Hence there can be at most one $4$-subtuple~$\vec{p}\,{}'$ for which the prof\/ile is unique.
So we can distinguish two cases:

\looseness=1
\emph{Case i.} For all $4$-subtuples~$\vec{p}\,{}'$, the prof\/ile is not uniquely determined.
Then by Lemma~\ref{lem:base_cases} we have that all points $\vec{p}\,{}'$ and $\vec{q}$ lie on a~cyclic cubic for each
$\vec{p}\,{}'$.
We can consider the linear system of all cyclic cubics through $q_1,\dotsc,q_m$.
It def\/ines a~rational map $g:(\mathbb P_{\mathbb C}^1 \times \mathbb P_{\mathbb C}^1)\dashrightarrow \mathbb P_{\mathbb
C}^N$, where~$N$ is the rank of the linear system (namely its dimension as a~projective space).
Since~$4$ real points in the plane give independent conditions for interpolating a~cyclic cubic, it follows that $N \le
3$.
Suppose now that the linear system giving~$g$ does not have any f\/ixed component.
Then because any four target points lie on cyclic cubic, it follows that any four of their images under~$g$ lie on
hyperplane.
Hence by a~linear algebra argument all images of the target points lie on a~hyperplane, and this hyperplane corresponds
to a~cyclic cubic curve as stated in the thesis.
If the linear system has some f\/ixed component then this f\/ixed component can only be a~circle or a~line, and this
contradicts the hypothesis, since we supposed that no~$4$ measure points are cocircular or aligned.

\emph{Case ii.} There exists a~$4$-subtuple~$\vec{p}\,{}'$ for which the prof\/ile is unique, but for any other
$4$-subtuple~$\vec{p}\,{}''$ the prof\/ile is not unique.
After possibly relabeling the points we can assume that $\vec{p}\,{}' = (p_2,p_3,p_4,p_5)$.
We can def\/ine a~map $g:(\mathbb P_{\mathbb C}^1 \times \mathbb P_{\mathbb C}^1)\dashrightarrow \mathbb P_{\mathbb C}^N$
as in Case~i.
Suppose that the linear system~$\Lambda$ of cyclic cubics through $q_1, \dotsc, q_m$ does not have any f\/ixed component.
Then by hypothesis we have that each $4$-tuple of points $\big(g(p_1), g(p_i), g(p_j), g(p_k) \big)$ lies on the same
plane for all $i,j,k \in \{2,3,4,5\}$.
This forces $g(p_2), \dotsc, g(p_5)$ to lie on the same plane, but this contradicts the hypothesis, since the prof\/ile
of $\vec{p}\,{}' = (p_2,p_3,p_4,p_5)$ is unique.
Hence the linear system~$\Lambda$ has a~f\/ixed component, but this means that the measure points are collinear or
cocircular.
Hence under our hypothesis Case~ii never arises.

If now~$t$ is arbitrary, we can argue as before showing that at most one $t-1$ tuple of target points admits a~unique
prof\/ile, and use induction on~$t$ since the hypothesis ensures that Case~ii never happens.

The proof for the dual assertion assuming non-uniqueness of the co-prof\/ile is similar, but with one modif\/ication.
If $m \ge 4$ and $t \ge 5$, then the linear system of cyclic cubics through all target points def\/ines the rational map
$g:(\mathbb P_{\mathbb C}^1 \times \mathbb P_{\mathbb C}^1) \dashrightarrow \mathbb P_{\mathbb C}^N$, where~$N$ is the
rank of the linear system (which, in this case, is less than or equal to~$2$).
Here is the modif\/ication: $5$ real points may fail to give independent conditions for interpolating a~cyclic cubic, and
this happens if and only if the all target points are collinear.
However this situation is ruled out by the hypothesis, so we have the thesis.
\end{proof}

\section{Point identif\/ication}
\label{sec:point_identification}

We are turning to the task of determining the target points $p_1, \dotsc, p_t$, assuming we already know the surface
$S_{\vec{p}}\subset \mathbb P_{\mathbb C}^{t-1}$, describing all possible results of measurements modulo~$\pi$.
This problem can only be solved up to similarities, namely rotations,
translations and dilations.
For $t = 3$, the problem is clearly not solvable, because in that case $S_{\vec{p}} = \mathbb P_{\mathbb C}^2$ does not
give any information on the target points.

\begin{Proposition}
\label{thm:t5pi}
If $t \ge 5$ and the target points are not collinear, then the point identification problem has a~unique solution.
\end{Proposition}
\begin{proof}
One sees that it is enough to prove the statement for $t = 5$.
To do this, one can proceed algebraically, showing that point identif\/ication boils down to solving systems of linear
equations admitting in general a~unique solution.
Otherwise one can use the following argument.
We assume that no 4 points are collinear.
If we consider four points $p_i$, $p_j$, $p_k$, $p_h$ out of the f\/ive, there exists a~unique measure point $q_{ijkh}$ such
that both angles $\measuredangle_{q_{ijkh}}(p_i, p_k)$ and $\measuredangle_{q_{ijkh}}(p_j, p_h)$ are zero
(modulo~$\pi$).
In fact this point is given by the intersection of the two diagonals $\overrightarrow{p_i p_k}$ and $\overrightarrow{p_j
p_h}$.
Then one performs the following construction: take the three tuples $(p_1, p_2, p_3, p_4)$, $(p_1, p_2, p_5, p_4)$ and
$(p_1, p_2, p_3, p_5)$ and denote by~$q_1$, $q_2$ and~$q_3$ the points $q_{1234}$, $q_{1254}$ and $q_{1235}$
respectively, as in Fig.~\ref{figure:point_identification}.

\begin{figure}[ht] \centering
\begin{overpic}
[width=0.25\textwidth]{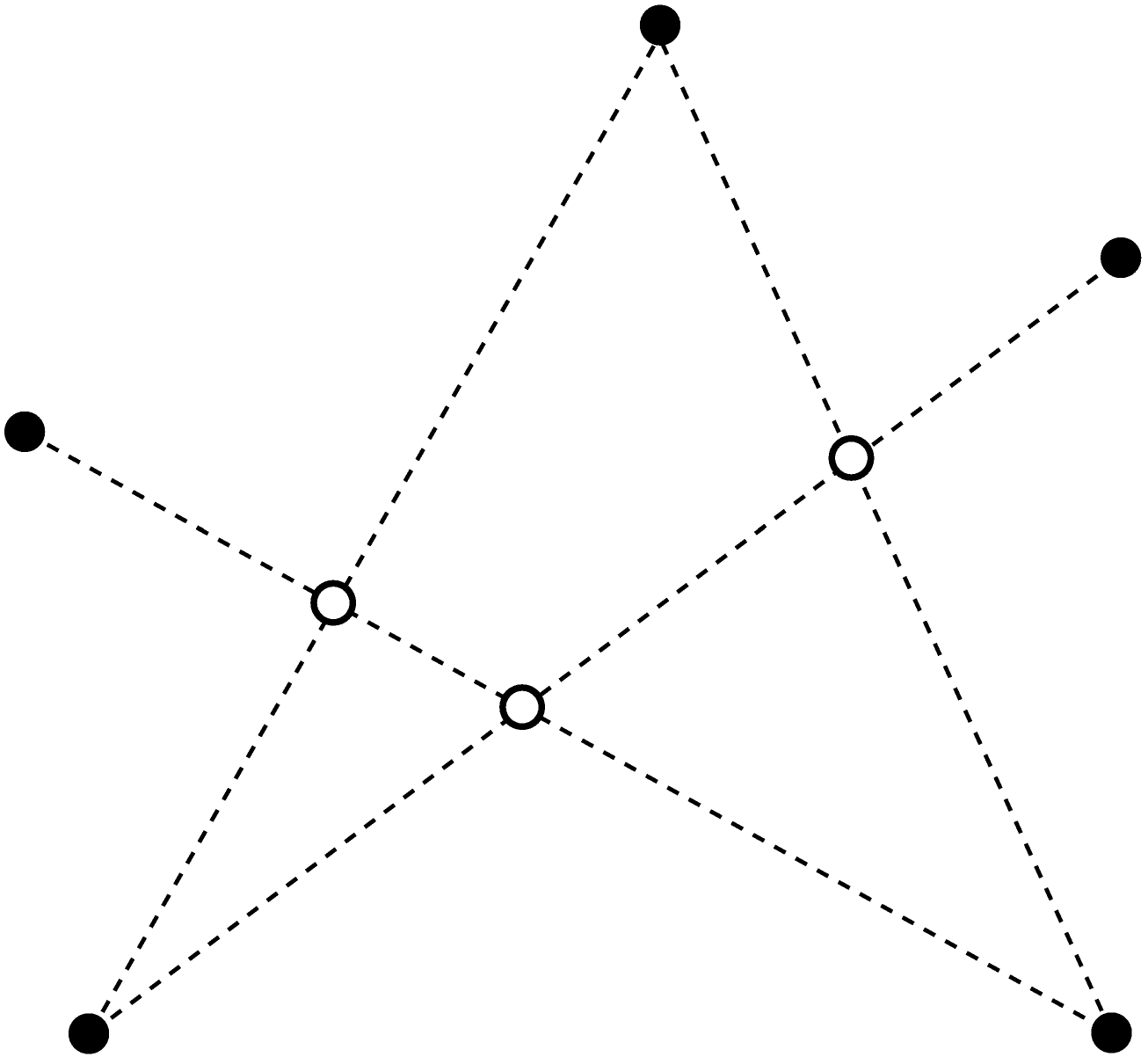}
\begin{small}
\put(-2,0){$p_1$} \put(102,0){$p_2$} \put(96,74){$p_3$} \put(-4,60){$p_4$} \put(46,88){$p_5$}

\put(25,45){$q_2$} \put(41,36){$q_1$} \put(75,59){$q_3$}

\put(27,30){$\beta$} \put(42,23){$\alpha$} \put(70,44){$\gamma$}
\end{small}
\end{overpic}
\caption{Identif\/ication of~$5$ points from the knowledge of their prof\/ile.}
\label{figure:point_identification}
\end{figure}
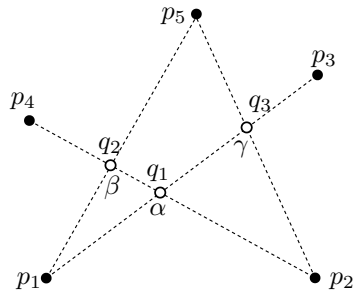

Since by hypothesis we know the prof\/ile of $p_1, \dotsc, p_5$, then we know the angles $\alpha =
\measuredangle_{q_1}(p_1, p_2)$, $\beta = \measuredangle_{q_2}(p_1, p_2)$ and $\gamma = \measuredangle_{q_3}(p_1, p_2)$
(modulo~$\pi$).
Then one sees that the angle $\measuredangle_{p_5}(p_1, p_2)$ is equal to $\beta + \gamma - \alpha$, namely it is
completely determined by the prof\/ile.
By symmetry, all angles $\measuredangle_{p_k}(p_i, p_j)$ are completely determined by the prof\/ile, hence the solution is
unique up to similarities.

If~$4$ points are collinear, and the~$5^{\mathrm{th}}$ is not, then the above construction can still be applied to
obtain enough angles, allowing to reconstruct the points up to similarity.
\end{proof}

\begin{Remark}
For collinear target points, the situation is also quite clear: two collinear tuples of measure points have the same
prof\/ile if and only if they are projectively equivalent.
\end{Remark}

\begin{Lemma}
\label{lem:m4pi}
If $m\ge 4$, then the point identification problem from the co-profile has a~unique solution.
\end{Lemma}

\begin{proof}
Let $m \ge 4$ and let $q_1, \dotsc, q_m$ be the measure points.
Let~$p_1$ be a~target point.
Assume that we know the co-prof\/ile surface $S'\subset\mathbb P_{\mathbb C}^{m}$ of $\vec{q} = (q_1,\dotsc,q_m)$ and
$p_1$.
Choose four points on $S'$ with af\/f\/ine coordinates on the complex unit circle, general with respect to this property.
Then the af\/f\/ine coordinates are the measurements of four additional target points $p_2$, $p_3$, $p_4$, $p_5$; by generality
assumptions, we may assume that there is no cyclic cubic curve passing through all target and measure points and that
the target points are neither collinear nor cocircular.
The transpose of the matrix consisting of the four coordinate vectors is now considered as an input for prof\/ile
interpolation.
Since $m = 4$, prof\/ile interpolation is uniquely solvable; let $S \subset \mathbb P_{\mathbb C}^3$ be the result.
By Proposition~\ref{thm:t5pi}, the prof\/ile determines $p_1,\dotsc,p_5$ uniquely (up to similarity).
The points $q_1,\dotsc,q_m$ can then be obtained by the resection map (see the proof of
Proposition~\ref{prop:measure_birational}).
\end{proof}

\begin{Proposition}
%\label{thm:t4pi}
If $t = 4$, then the point identification problem for the profile has, in general, two solutions.
\end{Proposition}

\begin{proof}
Assume that $p_1, \dotsc, p_4$ are unknown target point in general position with a~known prof\/ile surface $S \subset
\mathbb P_{\mathbb C}^3$.
As in the previous proof, we choose four points on~$S$ with af\/f\/ine coordinates on the complex unit circle, general with
respect to this property.
Their af\/f\/ine coordinates are the measurements from four unknown measurement points $q_1, \dotsc, q_4$.
We transpose the mat\-rix consisting of these four coordinate vectors and get an input for the co-prof\/ile interpolation
problem.
By Proposition~\ref{prop:m4int}, we get two solutions $S'$ and $S''$ in $\mathbb P_{\mathbb C}^3$.
For each of these two, using Lemma~\ref{lem:m4pi} we can then identify the points $\vec{q} = (q_1,\dotsc,q_m)$ and $p_1$
uniquely, up to similarity; by permuting target points, we can get also the remaining target points for both cases.
\end{proof}

\begin{notation}
For more f\/luent language, we call a~sequence of pairwise distinct $4$ points a~{\em quadrilateral}, and the second
solution to its point identif\/ication problem its {\em twin} (which is determined up to similarity).
\end{notation}

\begin{Remark}
\label{rem:twin}
Here is a~geometric construction for the twin quadrilateral.
If the vertices of a~quadrilateral~$V$ are cocircular, then~$V$ is its own twin (in other words, it can be recognized
uniquely by point identif\/ication).
Otherwise, the twin is obtained by constructing the centers of the circles that are circumscribed to the four triangles
formed by the vertices of~$V$ and then applying a ref\/lection to this quadrilateral (see Fig.~\ref{figure:twins}).

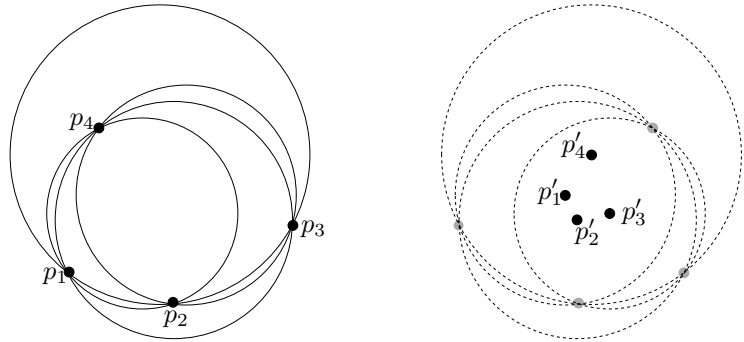
\begin{figure}[ht!] \centering
\begin{tabular}{ccc}
%\subfigure{
\begin{overpic}
[width=0.25\textwidth]{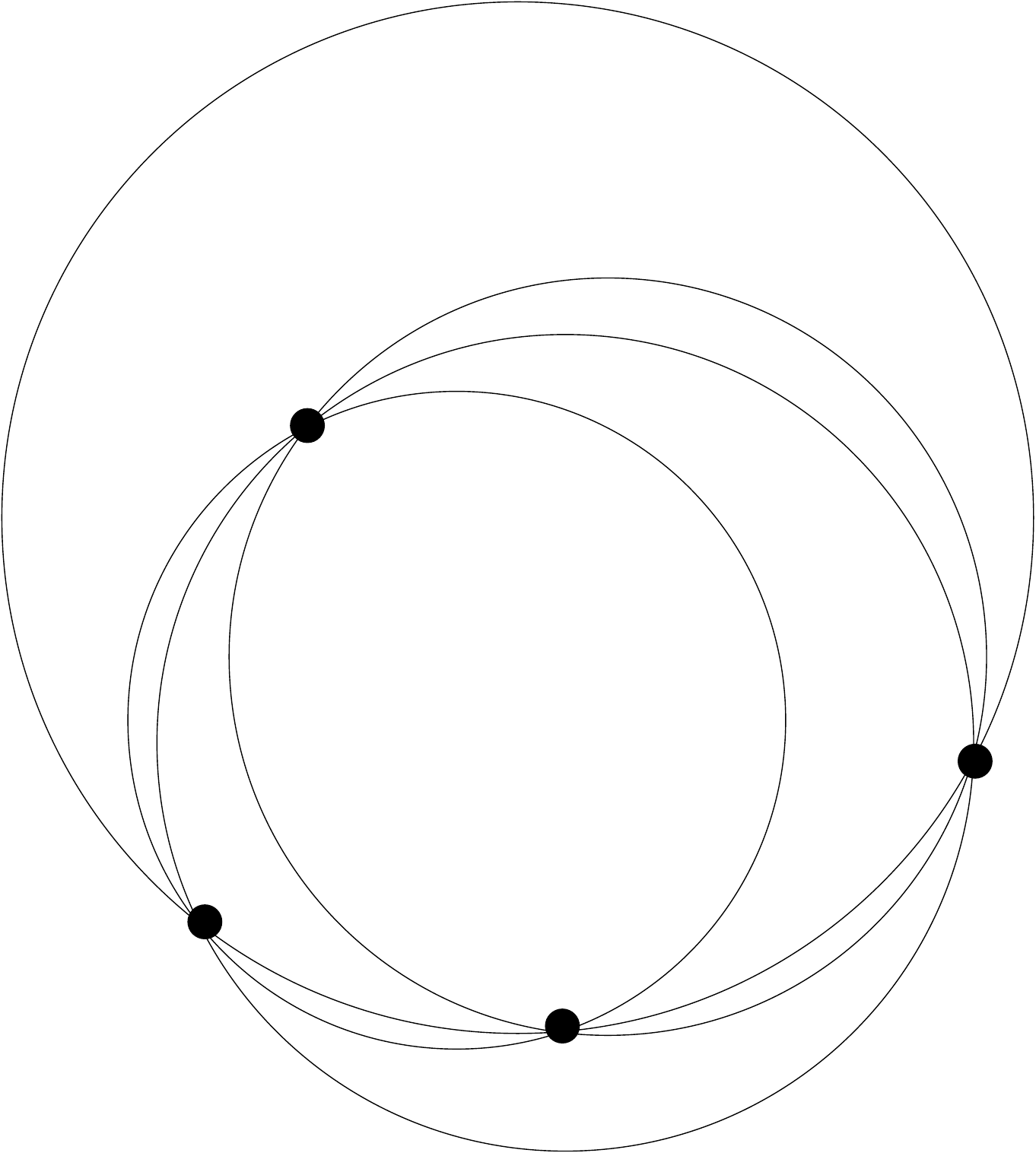}
\begin{small}
\put(10,18){$p_1$} \put(46,5){$p_2$} \put(87,32){$p_3$} \put(18,64){$p_4$}
\end{small}
\end{overpic}
%}
\hspace{0.4cm} & \hspace{0.4cm} & %\subfigure{
\begin{overpic}
[width=0.25\textwidth]{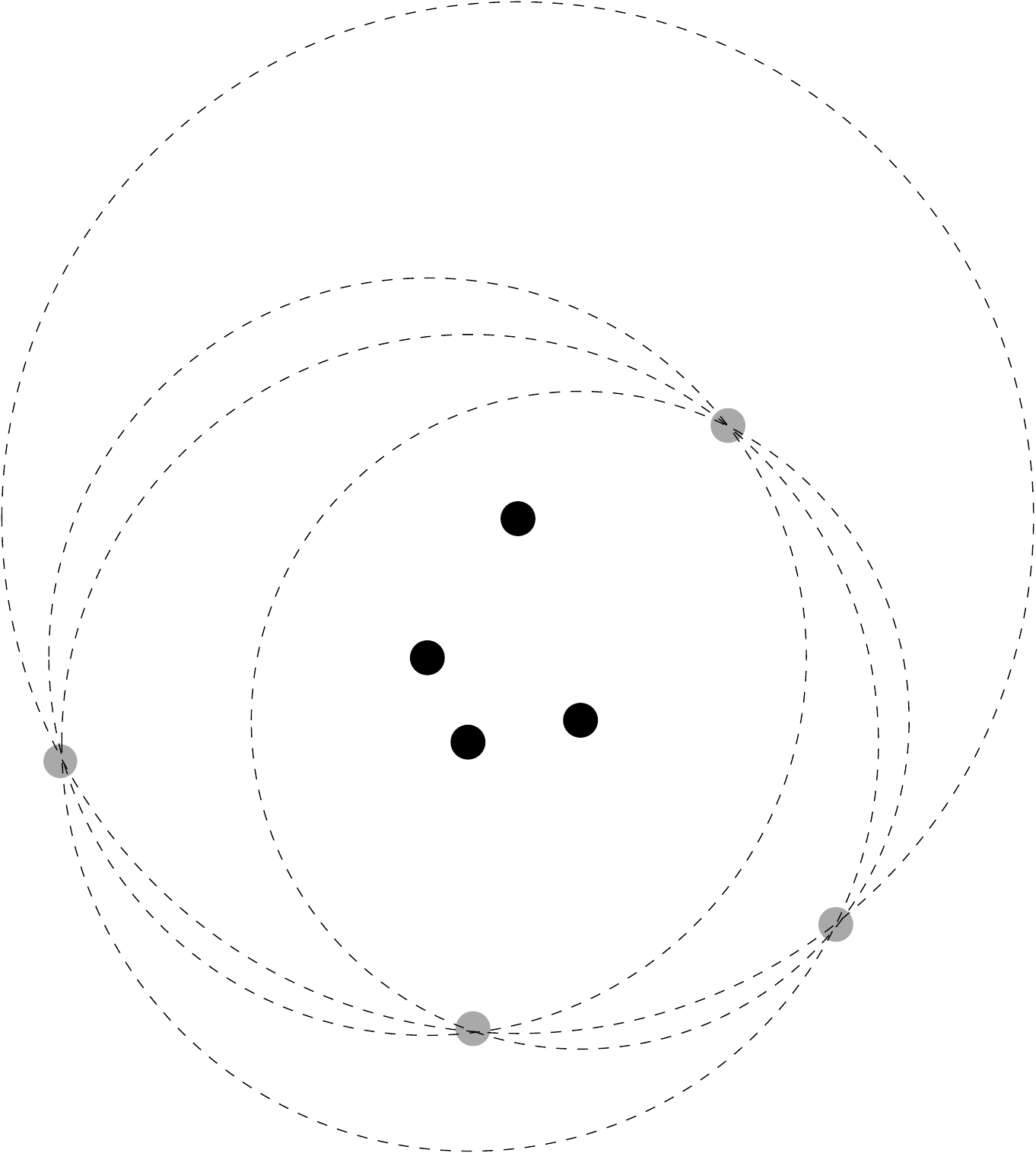}
\begin{small}
\put(54,36){$p_3'$} \put(40,30){$p_2'$} \put(36,56){$p_4'$} \put(29,42){$p_1'$}
\end{small}
\end{overpic}
%}
\end{tabular}
\caption{Construction of the twin quadrilateral: given $p_1$, $p_2$, $p_3$ and $p_4$ we consider the four circles passing
through 3 out of the 4 points.
The centers of these circles give the vertices $p_1'$, $p_2'$, $p_3'$ and $p_4'$ of the mirrored twin quadrilateral.}
\label{figure:twins}
\end{figure}

The proof of this statement is computational: we compute the prof\/ile for a~quadrilateral with symbolic coordinates, and
compare with the prof\/ile for the centers of the four circles above.
This was done using the computer algebra system Maple\footnote{See
\url{http://www.risc.jku.at/people/jschicho/pub/twinmaple.html}.}.
\end{Remark}

We summarize the results we obtained so far about point identif\/ication in the following theorem.

\begin{Theorem}
\label{thm:sup}
Assume that $(t-3)(m-2)\ge 2$.
The problem of identifying unknown target points $p_1,\dotsc,p_t$ and measure points $q_1,\dots,q_m$ from a~given double
angle matrix $M_{\vec{p},\vec{q}}$ is uniquely solvable unless we are in the following cases:
\begin{itemize}\itemsep=0pt
\item $t=4$ $($here there are, in general, two solutions$)$;
\item $m=3$ $($here there are, in general, two solutions$)$;
\item all target and measure points lie on a~cyclic cubic curve;
\item at least $5$ target points are cocircular or collinear;
\item at least $4$ measure points are cocircular or collinear.
\end{itemize}
\end{Theorem}

\section{Ambiguities for directed angles}
\label{sec:twinmap}

In this last section we give a~closer study of the case of~$4$ target points, and in particular we try to highlight some
properties relating a~quadrilateral to its twin.

The complex entries of the double angle matrix $M_{\vec{p},\vec{q}}$ give conditions on the lines through~$p_i$
and~$q_j$ for $i \in \{1, \dotsc, t\}$ and $j \in \{1, \dotsc, m\}$, but they do not include information about the
directions, i.e.\ on which of the two rays the point~$p_i$ is lying.
In a~situation where we have two solutions for the point identif\/ication problem, it may still be possible to tell the
right solution by taking direction information into account.
However our f\/inal result (see Theorem~\ref{thm:real}) will make explicit that if measurements are taken from points in
some prescribed region of the plane, then a~quadrilateral cannot be distinguished from its twin even if we take
directions into account.

\begin{Remark}
If the point identif\/ication has inf\/initely many solutions, then there are always ambiguities even if we take direction
information into account, because the direction information only allows to tell f\/initely many cases apart.
\end{Remark}

Let $\vec{p} = (p_1, \dotsc, p_4)$ be a~{\em nondegenerated quadrilateral}, i.e.\ a~quadruple of points such that no
three are collinear.
Then from what we saw at the end of Section~\ref{sec:point_identification} there is a~twin quadrilateral $\vec{p}\,{}' =
(p'_1,\dotsc,p'_4)$ and both $\vec{p}$ and $\vec{p}\,{}'$ have the same prof\/ile surface $S \subset \mathbb P_{\mathbb C}^3$.
We def\/ine the {\em twin map}:
\begin{gather*}
\rho_{\vec{p},\vec{p}\,{}'}:  \ \mathbb R^2 \setminus \{p_1,\dotsc,p_4\}  \longrightarrow  \mathbb R^2,
\\
\hphantom{\rho_{\vec{p},\vec{p}\,{}'}: \ }~q   \mapsto   \rho_{\vec{p},\vec{p}\,{}'}(q)  :=   \big(f_{\vec{p}\,{}'}^{-1} \circ f_{\vec{p}}\big)(q).
\end{gather*}
If $p_1, \dotsc, p_4$ are cocircular, then the twin quadrilateral is similar to the original one and the twin map is
just a~composition of a~rotation, a~dilation and a~translation.
From now on, we assume that the points are not cocircular.

Algebraically, the twin map can be extended to a~birational automorphism of~$\mathbb P^1_{\mathbb C} \times \mathbb
P^1_{\mathbb C}$.
The exceptional curves partition the real plane minus the exceptional locus into open connected regions.
Each open region is mapped homeomorphically to its image region.
These image regions are the result of an analogous partition by the exceptional curves of~$\rho_{\vec{p}\,{}'\!, \vec{p}}$.
This strongly motivates our interest in the exceptional curves of the twin map.
They are described in the next proposition.

\begin{Proposition}
\label{prop:exc}
The exceptional curves of~$\rho_{\vec{p},\vec{p}\,{}'}$ are the four circles passing through three of the points $p_1,
\dotsc, p_4$.
\end{Proposition}

\begin{proof}
The base locus of the twin map is equal to the set $\{p_1,\dotsc,p_4\}$.
Because the image is a~smooth surface, the exceptional curves are in ($-1$)-classes of the blowing up~$\mathbb{Y}$
of~$\mathbb P^1_{\mathbb C}\times\mathbb P^1_{\mathbb C}$ at $p_1, \dotsc, p_4$.
The class group of~$\mathbb{Y}$ is freely generated by the class~$L_1$ of bidegree $(1,0)$, the class~$L_2$ of bidegree
$(0,1)$, and the four exceptional classes $E_1$, $E_2$, $E_3$, $E_4$ of the blowing up map.
Here is the list of all ($-1$)-classes in the class group of~$\mathbb{Y}$:
\begin{itemize}\itemsep=0pt
\item  $E_1,\dotsc,E_4$, corresponding to the four base points;
\item  $L_1-E_1,\dotsc,L_1-E_4,L_2-E_1,\dotsc,L_2-E_4$, corresponding to the f\/ibers of the two projections
passing through the base points;
\item  $L_1+L_2-E_2-E_3-E_4,\dotsc,L_1+L_2-E_1-E_2-E_3$, corresponding to the circles passing through three of
the four base points.
\end{itemize}
This list is very well known, since the blow up of~$\mathbb P^1_{\mathbb C}\times\mathbb P^1_{\mathbb C}$ at four points
is isomorphic to the blow up of $\mathbb P^2_{\mathbb C}$ at f\/ive points.
Hence $\mathbb{Y}$ is a~Del Pezzo surface of degree~$4$, and it contains exactly the $16$ $(-1)$-classes listed before
(for a~reference, see~\cite[Theorems~24.3, 24.4 and~24.5]{Manin} or~\cite[Chapter~8 and in particular the
remark following Lemma~8.2.22]{Dolgachev}).
In the factorization of $\rho_{\vec{p},\vec{p}\,{}'}$ into blowing ups and blowing downs, we have four blowing ups, so we
also need to have four blowing downs.
The only choice of four ($-1$)-classes not containing the exceptional classes of the blowing up, and such that it is symmetric
under permutations of $E_1, \dotsc, E_4$ and under permutations of~$L_1$ and~$L_2$, is the choice of the four circles.
\end{proof}

We denote by $C_1, \dotsc, C_4$ the four circles mentioned in Proposition~\ref{prop:exc}, where~$C_i$ is the circle that
does not pass through~$p_i$, for $i \in \{1, \dotsc, 4\}$.
We call them the {\em fundamental circles} of~$\vec{p}$.

The following proposition is not really needed; we just mention it because it is a~nice description of the structure of
the partition by fundamental circles.

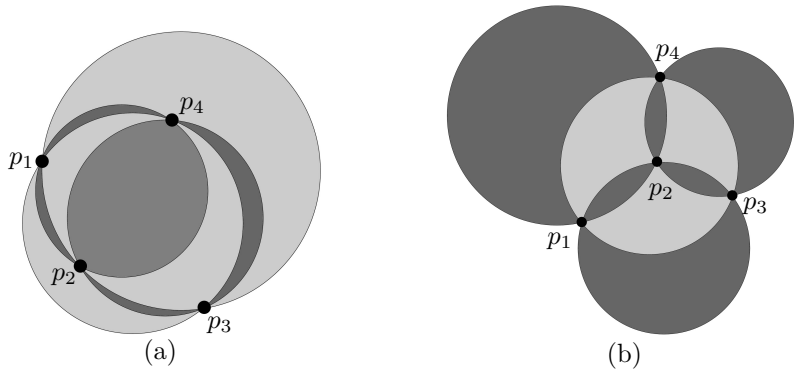
\begin{figure}[ht!] \centering
\begin{tabular}{cc} %\subfigure[][]{
\begin{overpic}
[width=0.25\textwidth]{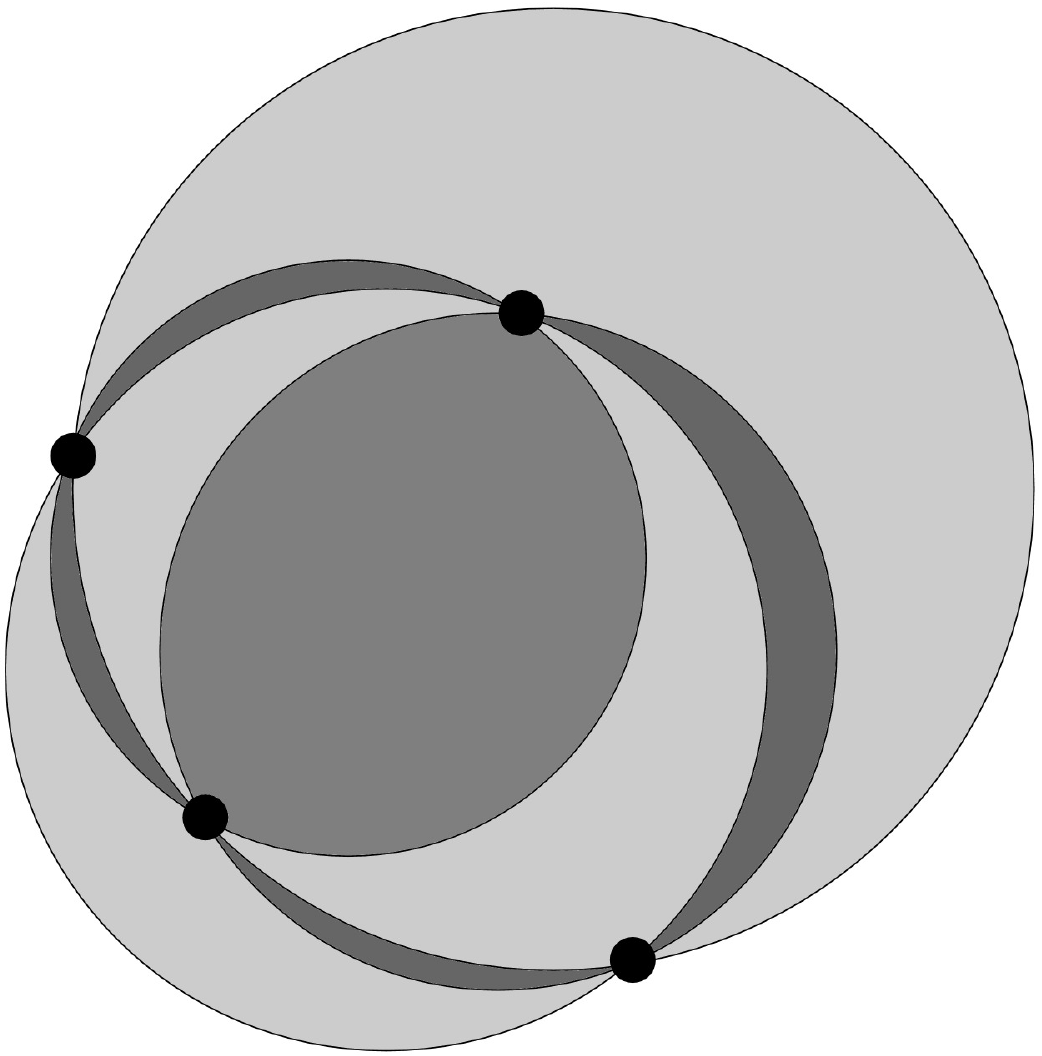}
\begin{small}
\put(-4,56){$p_1$} \put(10,18){$p_2$} \put(61,2){$p_3$} \put(52,74){$p_4$} \put(40,-8){(a)}
\end{small}
\end{overpic}
%\label{figure:regions_convex}
%}
\hspace{0.5cm} & \hspace{0.5cm} %\subfigure[][]{
\begin{overpic}
[width=0.29\textwidth]{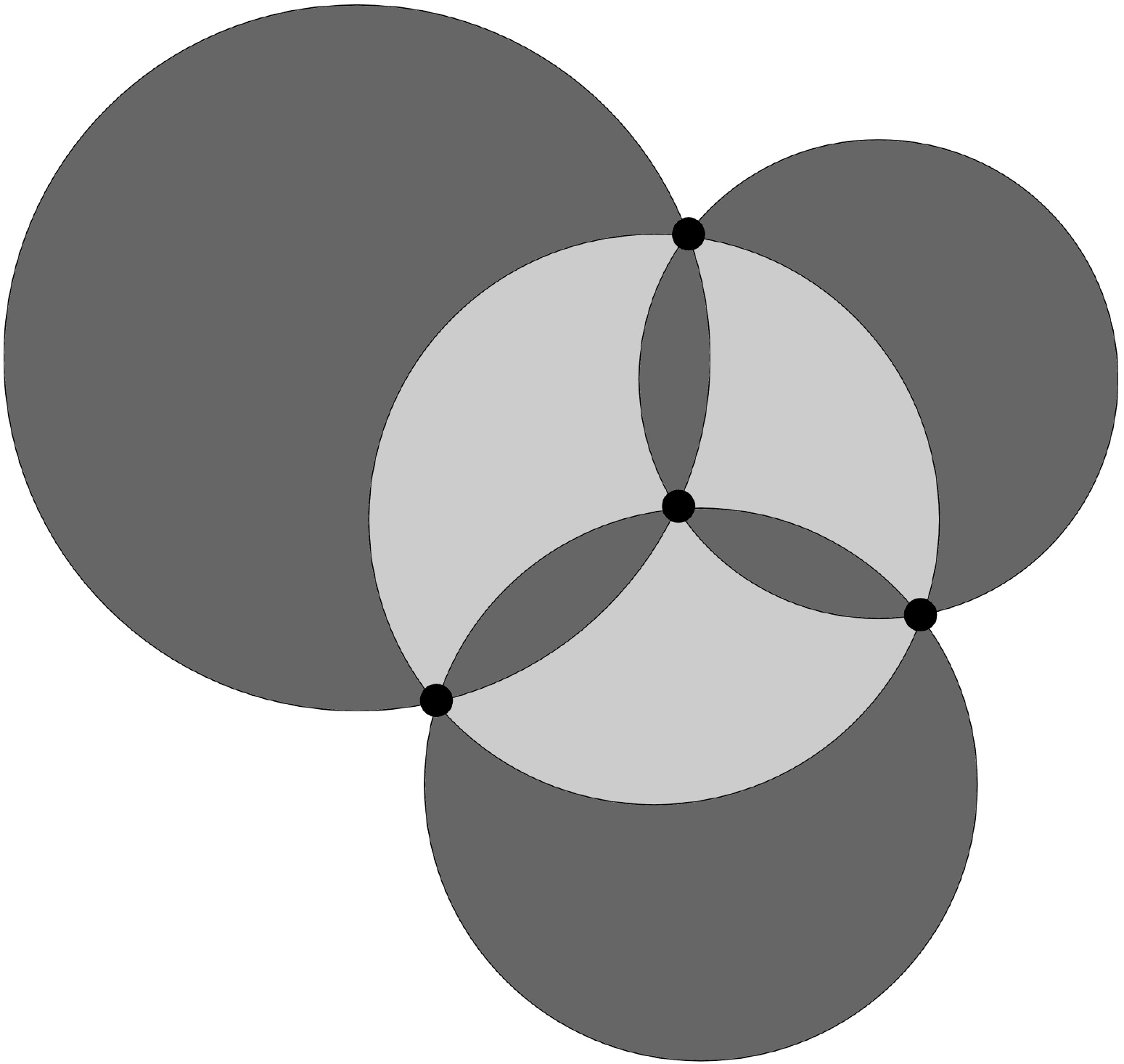}
\begin{small}
\put(29,26){$p_1$} \put(58,40){$p_2$} \put(85,36){$p_3$} \put(60,80){$p_4$} \put(46,-8){(b)}
\end{small}
\end{overpic}
%\label{figure:regions_concave}
%}
\end{tabular}\vspace{1mm}
\caption{The~$10$ regions forming the complement of the fundamental circles of a~quadrilateral: in the convex
case~(a) the boundary of the unbounded region is given by $2$ arcs, while in the concave
case~(b) the boundary is given by~$3$ arcs.}
\label{figure:regions}
\end{figure}

\begin{Proposition}
The complement of the fundamental circles of a~quadrilateral~$\vec{p}$ has~$10$ connected components.
Four of them are delimited by three circles, and the other~$6$ are delimited by two circles.
The boundary of the unique unbounded region is constituted by~$2$ arcs if and only if the quadrilateral is convex $($see
Fig.~{\rm \ref{figure:regions})}.
\end{Proposition}

If $\vec{p}$ is convex, then the unique region containing the intersection of the inner diagonals is called the {\em
inner region}.
It is important to observe that the twin of a~convex quadrilateral is again convex, and the twin map maps the inner
region to the inner region.

\begin{Theorem}
\label{thm:real}
Let~$\vec{p}$ be a~nondegenerate convex quadrilateral and assume that all measure points $q_1,\dotsc,q_m$ $($where~$m$ is
arbitrary$)$ are contained in the inner region.
Let $\vec{p}\,{}'$ be the twin quadrilateral of~$\vec{p}$ and set $q'_i = \rho_{\vec{p},\vec{p}\,{}'}(q_i)$ for $i \in \{1,
\dotsc, m\}$.
Then we have
\begin{gather*}
\measuredangle_{q_j}(p_1,p_{i})
=
\measuredangle_{q'_j}(p'_1,p'_{i})
\end{gather*}
for all $i \in \{2, \dotsc,t\}$ and $j \in \{1, \dotsc, m\}$.
\end{Theorem}

\begin{proof}
The construction guarantees that $M_{\vec{p},\vec{q}} = M_{\vec{p}\,{}',\vec{q}\,{}'}$.
Hence for all $i \in \{2, \dotsc, t\}$ and $j \in \{1, \dotsc, m\}$ we have the equality
$\big(\measuredangle_{q_j}(p_1,p_{i})\big)^2 = \big(\measuredangle_{q'_j}(p'_1,p'_{i})\big)^2$, where we represent
angles by complex numbers, so that doubling the angle corresponds to squaring the complex representation.
Hence the quotient of left hand side and right hand side is $\pm 1$.
We consider this quotient as a~function~$\phi$ from the complement of the fundamental circles to~$\{1,-1\}$.
It is easy to see that this map is continuous, therefore constant on any region, since they are connected.
At the intersection of the inner diagonals the map~$\phi$ assumes the value~$1$, therefore it is constantly equal to~$1$
on the whole inner region.
From this the thesis follows.
\end{proof}

\subsection*{Acknowledgements}

We thank Kristian Ranestad for pointing out to us the solution of a~crucial interpolation problem, Bert J\"uttler and
Martin Peternell for helping us to orient ourselves in the literature on geodesy, and the anonymous referees for helpful
comments and pointers to literature on photogrammetry.
The authors' research is supported by the Austrian Science Fund (FWF): W1214-N15/DK9 and P26607~-- ``Algebraic Methods
in Kinematics: Motion Factorisation and Bond Theory''.

\pdfbookmark[1]{References}{ref}
\LastPageEnding

\end{document}